\magnification=1200

%
%
\newdimen\FigSize	\FigSize=.9\hsize 
%
\newskip\abovefigskip	\newskip\belowfigskip
\gdef\epsfig#1;#2;{\par\vskip\abovefigskip\penalty -500
   {\everypar={}\epsfxsize=#1\nd
    \centerline{\epsfbox{#2}}}%
    \vskip\belowfigskip}%
%
\newskip\figtitleskip
\gdef\tepsfig#1;#2;#3{\par\vskip\abovefigskip\penalty -500
   {\everypar={}\epsfxsize=#1\nd
    \vbox
      {\centerline{\epsfbox{#2}}\vskip\figtitleskip
       \centerline{\figtitlefont#3}}}%
    \vskip\belowfigskip}%
%
\newcount\FigNr	\global\FigNr=0
\gdef\nepsfig#1;#2;#3{\global\advance\FigNr by 1
   \tepsfig#1;#2;{Figure\space\the\FigNr.\space#3}}%
%
%
%
\gdef\ipsfig#1;#2;{
   \midinsert{\everypar={}\epsfxsize=#1\nd
	      \centerline{\epsfbox{#2}}}%
   \endinsert}%
%
\gdef\tipsfig#1;#2;#3{\midinsert
   {\everypar={}\epsfxsize=#1\nd
    \vbox{\centerline{\epsfbox{#2}}%
          \vskip\figtitleskip
          \centerline{\figtitlefont#3}}}\endinsert}%
%
\gdef\nipsfig#1;#2;#3{\global\advance\FigNr by1%
  \tipsfig#1;#2;{Figure\space\the\FigNr.\space#3}}%
\newread\epsffilein    
\newif\ifepsffileok    
\newif\ifepsfbbfound   
\newif\ifepsfverbose   
\newdimen\epsfxsize    
\newdimen\epsfysize    
\newdimen\epsftsize    
\newdimen\epsfrsize    
\newdimen\epsftmp      
\newdimen\pspoints     
\pspoints=1bp          
\epsfxsize=0pt         
\epsfysize=0pt         
\def\epsfbox#1{\global\def\epsfllx{72}\global\def\epsflly{72}%
   \global\def\epsfurx{540}\global\def\epsfury{720}%
   \def\lbracket{[}\def\testit{#1}\ifx\testit\lbracket
   \let\next=\epsfgetlitbb\else\let\next=\epsfnormal\fi\next{#1}}%
\def\epsfgetlitbb#1#2 #3 #4 #5]#6{\epsfgrab #2 #3 #4 #5 .\\%
   \epsfsetgraph{#6}}%
\def\epsfnormal#1{\epsfgetbb{#1}\epsfsetgraph{#1}}%
\def\epsfgetbb#1{%
%
%
\openin\epsffilein=#1
\ifeof\epsffilein\errmessage{I couldn't open #1, will ignore it}\else
%
%
   {\epsffileoktrue \chardef\other=12
    \def\do##1{\catcode`##1=\other}\dospecials \catcode`\ =10
    \loop
       \read\epsffilein to \epsffileline
       \ifeof\epsffilein\epsffileokfalse\else
%
%
          \expandafter\epsfaux\epsffileline:. \\%
       \fi
   \ifepsffileok\repeat
   \ifepsfbbfound\else
    \ifepsfverbose\message{No bounding box comment in #1; using
defaults}\fi\fi
   }\closein\epsffilein\fi}%
%
%
\def\epsfsetgraph#1{%
   \epsfrsize=\epsfury\pspoints
   \advance\epsfrsize by-\epsflly\pspoints
   \epsftsize=\epsfurx\pspoints
   \advance\epsftsize by-\epsfllx\pspoints
%
%
   \epsfxsize\epsfsize\epsftsize\epsfrsize
   \ifnum\epsfxsize=0 \ifnum\epsfysize=0
      \epsfxsize=\epsftsize \epsfysize=\epsfrsize
%
arithmetic!
%
     \else\epsftmp=\epsftsize \divide\epsftmp\epsfrsize
       \epsfxsize=\epsfysize \multiply\epsfxsize\epsftmp
       \multiply\epsftmp\epsfrsize \advance\epsftsize-\epsftmp
       \epsftmp=\epsfysize
       \loop \advance\epsftsize\epsftsize \divide\epsftmp 2
       \ifnum\epsftmp>0
          \ifnum\epsftsize<\epsfrsize\else
             \advance\epsftsize-\epsfrsize \advance\epsfxsize\epsftmp
\fi
       \repeat
     \fi
   \else\epsftmp=\epsfrsize \divide\epsftmp\epsftsize
     \epsfysize=\epsfxsize \multiply\epsfysize\epsftmp
     \multiply\epsftmp\epsftsize \advance\epsfrsize-\epsftmp
     \epsftmp=\epsfxsize
     \loop \advance\epsfrsize\epsfrsize \divide\epsftmp 2
     \ifnum\epsftmp>0
        \ifnum\epsfrsize<\epsftsize\else
           \advance\epsfrsize-\epsftsize \advance\epsfysize\epsftmp \fi
     \repeat
   \fi
%
%
   \ifepsfverbose\message{#1: width=\the\epsfxsize,
height=\the\epsfysize}\fi
   \epsftmp=10\epsfxsize \divide\epsftmp\pspoints
   \vbox to\epsfysize{\vfil\hbox to\epsfxsize{%
      \includegraphics{#1}%
      \hfil}}%
\epsfxsize=0pt\epsfysize=0pt}%
%
%
{\catcode`\%=12
\global\let\epsfpercent=
%
%
\long\def\epsfaux#1#2:#3\\{\ifx#1\epsfpercent
   \def\testit{#2}\ifx\testit\epsfbblit
      \epsfgrab #3 . . . \\%
      \epsffileokfalse
      \global\epsfbbfoundtrue
   \fi\else\ifx#1\par\else\epsffileokfalse\fi\fi}%
%
%
\def\epsfgrab #1 #2 #3 #4 #5\\{%
   \global\def\epsfllx{#1}\ifx\epsfllx\empty
      \epsfgrab #2 #3 #4 #5 .\\\else
   \global\def\epsflly{#2}%
   \global\def\epsfurx{#3}\global\def\epsfury{#4}\fi}%
%
%
\def\epsfsize#1#2{\epsfxsize}%
%
%

\epsfverbosetrue			
\abovefigskip=\baselineskip		
\belowfigskip=\baselineskip		
\global\let\figtitlefont\bf		
\global\figtitleskip=.5\baselineskip	

\font\tenmsb=msbm10   
\font\sevenmsb=msbm7
\font\fivemsb=msbm5
\newfam\msbfam
\textfont\msbfam=\tenmsb
\scriptfont\msbfam=\sevenmsb
\scriptscriptfont\msbfam=\fivemsb
\def\Bbb#1{\fam\msbfam\relax#1}
\def\gtrsim{\mathop{\lower.5ex\hbox{$\buildrel
>\over{\scriptstyle\sim}$}}\nolimits}
\let\nd\noindent 
\def\qed{\hbox{\hskip 6pt\vrule width6pt height7pt depth1pt \hskip1pt}}
\def\natural{{\rm I\kern-.18em N}}

\def\integer{{\rm Z\kern-.32em Z}}
\def\chix{{\raise.5ex\hbox{$\chi$}}}
\def\Z{{\Bbb Z}}
\def\real{{\rm I\kern-.2em R}}
\def\R{{\Bbb R}}
\def\complex{\kern.1em{\raise.47ex\hbox{
            $\scriptscriptstyle |$}}\kern-.40em{\rm C}}

\def\H{{\Bbb H}}
\def\L{{k}} 
\def\E{{\Bbb E}}
\def\S{{\Bbb S}}
\def\G{{\cal G}}
\def\M{{\cal M}}
\def\O{{\cal O}}
\def\B{{\cal B}} 
\def\P{{\cal P}} 
\def\W{{W}} 
\def\T{{\cal T}} 
\def\t{{\Theta}} 
\def\vs#1 {\vskip#1truein}
\def\hs#1 {\hskip#1truein}

\def\Month{\ifcase\number\month \relax\or January \or February \or
  March \or April \or May \or June \or July \or August \or September
  \or October \or November \or December \else \relax\fi }
\def\date{\Month \the\day, \the\year}

\font\univfont=cmr8 scaled 1440
\font\deptfont=cmr6 scaled 1440
\font\addressfont=cmsl8
\def\utletterhead
{\hbox to \hsize{\parskip=4pt \parindent=0pt
   \vbox{\hbox to .1\hsize {\hss\sealfont\char0}}\hskip.05\hsize
   \vbox{\deptfont Department of Mathematics \endgraf
         \univfont The University of Texas \endgraf
         \vskip 6pt\hrule width .85\hsize \vskip 2pt
         \addressfont Austin, Texas 78712 $\cdot$ (512) 471-7711
                      $\cdot$ FAX (512) 471-9038
         \vskip 8pt}}}

  \hsize=6truein        \hoffset=.25truein 
  \vsize=8.8truein      
  \pageno=1     \baselineskip=12pt
  \parskip=0 pt         \parindent=20pt
  \overfullrule=0pt     \lineskip=0pt   \lineskiplimit=0pt
  \hbadness=10000 \vbadness=10000 
\nd
     \pageno=0
     
     \footline{\ifnum\pageno=0\hss\else\hss\tenrm\folio\hss\fi}
     \hbox{}
     \vskip 1truein\centerline{{\bf OPTIMALLY DENSE PACKINGS OF}}
     \vskip .1truein\centerline{{\bf HYPERBOLIC SPACE}}
     \vskip .2truein\centerline{by}
     \vskip .2truein
\centerline{{Lewis Bowen$^\dagger$
\footnote{*}{Research supported in part by NSF VIGRE grant No. DMS-0091946}}
\ \ and\ \  {Charles Radin}
\footnote{**}{Research supported in part by NSF Grant No. DMS-0071643 and
\vs0 Texas ARP Grant 003658-158\hfil}}

     \vskip .5truein\centerline{
     Mathematics Department, University of Texas at Austin}
     \vs.2
$^\dagger$ Current address: Mathematics Department, University of California
at Davis
\vs.5
     \centerline{{\bf Abstract}}
     \vs.1 \nd
In previous work a probabilistic approach to
controlling difficulties of density in hyperbolic space led to a
workable notion of optimal density for packings of bodies. In this
paper we extend an ergodic theorem of Nevo to provide an appropriate
definition of those packings to be considered optimally dense.
Examples are given to
illustrate various aspects of the density problem, in particular
the shift in emphasis from the analysis of individual packings to
spaces of packings.
     \vs.8
     \vs.2
     \centerline{Subject Classification:\ \ 52A40, 52C26, 52C23}
     \vfill\eject
\nd
{\bf 0. Introduction}
\vs.1
One of the main themes in discrete geometry has always been the study
of optimally dense packings of bodies, in regions of finite and,
especially, infinite volume -- the latter being the situation with
which we will be mainly concerned. There is a large literature for
packing in Euclidean spaces $\E^n$, but much less for packing in
hyperbolic spaces $\H^n$; see [GrW] and references therein, as well as
the classics [Fe5, Rog]. The difficulties in hyperbolic spaces have
been documented in many papers dating back at least to the early 50's
(see [Fe1-5], [Bo1-2], [BoF], [FeK], [FKK], [Kup]), and are well
understood to be related to the phenomenon of the exponential rate of
growth of the volume of a sphere with respect to radius.

In general it is difficult, even in $\E^n$, to actually determine the
optimal packing density even for simple shapes; for instance this is
unknown even for spheres for dimension $n\ge 4$. For certain polyhedra
such optima are occasionally computable, and one outgrowth of this,
since the late 60's, has been the subject sometimes referred to as
``aperiodic tiling''; see [Rad] and references therein. That subject
is concerned with the geometric symmetries of those packings (tilings)
which achieve optimal density, symmetries of unusual form, which are
studied through probabilistic techniques.

In hyperbolic space the difficulties in analyzing optimal density are,
as hinted above, much worse; not only is it hard to determine the
optimal density, say for spheres, it has even been hard to decide
whether this is a meaningful question.

Finding a meaningful, and computationally useful, way to analyze
optimal density in $\H^n$ was the subject of [BoR, Bow], in which some
of the probabilistic methods of aperiodic tiling were employed. The
ergodic theory approach to density in [BoR], summarized below, is a
modification of an approach based on Mass Transport suggested to us by
Oded Schramm (see [BeS]). We found, using ergodic theorems of Nevo et
al, that the quantity he suggested as a ``density'' is
in fact, in a statistical sense discussed below, the highest true
density - the limiting fraction of space covered by bodies in a fixed
packing. The papers [BoR, Bow] did not, however, convincingly address
the question of optimally dense packings themselves; that is, a
formalism for analyzing the density was found, but not a way to
address the packings that achieve that optimal density. That is one of
the main aims of this paper.

The other major objective of this paper is to present a series of
examples, some new and some which are variants of historically
important examples, adjusted to reveal significant connections.
Although one source of difficulty in the subject is well understood,
the rate of growth of volume with diameter, these examples point to a
very different source, related to structures in the space of packings.

\vs.2 \nd
{\bf I. Some troublesome packings}
\vs.1
We will be analyzing the density of certain subsets of Euclidean
$n$-dimensional space $\E^n$ or hyperbolic $n$-dimensional space $\H^n$
of
curvature $-1$; we let $\S$ stand for any of these spaces. The subsets
of $\S$ whose density we will consider will generally be ``packings''
of (infinitely many) ``bodies'' $\beta_j$, where a packing is a
collection
of bodies with pairwise disjoint interiors, and a body is a compact,
connected set which is the closure of its interior. One of the
features of our analysis will be an emphasis on distinguishing between
the density of the packing versus the density of the set which is the
union of the bodies in the packing; that is, it will be significant to
maintain the individuality of each of those bodies.

Our main focus will be on the ``densest'' packings possible by the
given bodies, and this requires examination of the primitive notion of
density. If we were packing a region $S$ of {\it finite} volume by the
bodies $\beta_j$, the density of such a packing would be unambiguous --
the fraction of the volume of $S$ covered by the bodies -- but density
must be defined more subtly for packings of a region, such as $\S$, of
infinite volume. The most widely accepted [FeK] primitive notion is
that the density of a packing $\P$ of $\S$ should be obtainable by choosing
a
family of finite volume regions $S_k$, with $S_k\subset S_{k+1}$ and
$\cup_k S_k=\S$, and the density of $\P$ should be

$$\lim_{k\to \infty} {vol(\P\cap S_k)\over vol(S_k)},\eqno(1)$$

\nd where $vol(\cdot )$ denotes volume in $\S$ and $\P\cap S_k$ denotes
the
portion of $S_k$ covered by bodies in $\P$. We would want the density
to be reasonably independent of the family $S_k$.

It is worth noting
that the limit in (1) can easily fail to exist. Consider the sequence
$\{D_j\,|\,j\ge 1\}$ of closed disks in $\E^2$, $D_j$ of radius $2^j$
and centered at the origin. Let $P_j$ be the annulus $D_j/D_{j-1}$
between successive disks, and let $S$ be the union of those $P_j$ with
$j\ge 2$ even. If we try to define the density of $S$ using the
expanding regions $S_k=D_k$, the sequence of local or approximate
densities $\displaystyle {vol(S\cap S_k)/vol(S_k)}$ would not
have a limit as $k\to\infty$, due to oscillation. (We could get the
same qualitative result by replacing our region $S$ by its
intersection with some simple packing of disks, such as the packing of
unit diameter disks whose centers have integer coordinates.)

Even though there are packings without a well defined density there is
no real difficulty in defining {\it optimal} density of packings in
Euclidean space. In fact we now show how to construct densest packings
of Euclidean space. Let $\S=\E^n$ and let $S_\L$ be a cube centered at
the origin, with edges of length $\L$ aligned with the axes. For any
$\L>0$, let $\P_\L$ be a packing by (congruent copies of the bodies
in) $\B\equiv \{\beta_j\}$ such that all bodies in $\P_\L$ intersect
$S_\L$ and $vol(\P_\L\cap S_\L)$ is optimally large.  (Such a packing
is easily shown to exist by a simple compactness argument [GrS;
p. 154].) For any packing $\P$ in $\S$ we define

$$d_\L(\P)={vol(S_\L\cap \P)\over vol(S_\L)}\eqno(2)$$

$$d_\L=\max_{\P} d_\L(\P)\eqno(3)$$

$$d=\limsup_{\L \to \infty} d_\L.\eqno(4)$$

\nd ($\lim_{\L \to \infty} d_\L$ exists but we do not need this fact.)

At this point it is convenient to have a space $\tilde \Sigma_\B$ of all
possible packings of $\S$ by the bodies $\beta_j$, equipped with a
metric topology such that a sequence of packings converges if and only if it
converges uniformly on compact subsets of $\S$. We will spell this out
in section II b, but assume for now such a space makes sense and is in
fact compact. Then we let $\P_\infty$ be an accumulation point of
$\{\P_\L\}$.
\vs.1
The following is a simple observation.
\vs.1 \nd
{\bf Lemma 1}.
$d_\L(g\P_\infty)\to d$ as $\L\to \infty$ for every fixed rigid motion $g$.
\vs.1 \nd
Proof. The main estimates needed are the simple facts, for $\L'> \L$:

$$vol(\P_\L\cap S_\L)\ge vol(\P_{\L'}\cap S_\L)\eqno(5)$$

$$vol(\P_{\L'}\cap S_\L)\ge vol(\P_\L\cap S_\L) -
[\L^n-(\L - C)^n],\eqno(6)$$

\nd where $C$ is larger than the diameter of any body in $\B$.
Equation (6) holds because if it did not one could arrive at a contradiction
by
altering $\P_{\L'}$ as follows.  First replace
the bodies of $\P_{\L'}$ that are completely contained in $S_\L$ by the
bodies of ${\P}_{\L}$ that do not overlap
the other bodies of ${P}_{\L'}$ (i.e.
that do not overlap any body of ${P}_{\L'}$ that overlaps the
complement of $S_\L$). Note that the volume of bodies of ${P}_\L$
that we have introduced is at least as large as the right hand side of
(6). Since $vol(\P_{\L'})$ is as large as possible, this operation
could not have increased its volume. This proves (6). Since
$\P_\infty$ is a limit of $\P_\L$ we get that (6) holds if $\P_{\L'}$ is
replaced by $\P_\infty$.

Finally, if $\L_m$ is a sequence such that $d_{\L_m}\to d$ as $m\to \infty$:

$$\eqalign{
|d-d_{\L_m}(g\P_\infty)|&\le |d_{\L_m} - d| + |d_{\L_m} -
d_{\L_m}(g\P_\infty)| \cr
&= |d_{\L_m} - d| + |d_{\L_m}(\P_{\L_m}) - d_{\L_m}(g\P_\infty)|}\eqno(7)$$

\nd and $|d_{\L_m}(\P_{\L_m}) - d_{\L_m}(g\P_\infty)|\to 0$ as $m\to \infty$
from (6).\qed
\vs.1
Thus, in Euclidean space optimally dense packings $\P$
exist for any collection $\B$ in the sense that their
density defined by (1) exists, and is as large as that for any packing.

As we shall see, the above technique does not extend to $\S=\H^n$
and therefore some other method must be used to define optimal density in
$\H^n$. Before exhibiting such a method, we present some examples to
highlight
some differences between hyperbolic and Euclidean packings.
\vs.1 \nd
{\bf Example 1 (half-space).}
\vs.1
Consider the half space region $S$, defined, in the upper half plane
model of the hyperbolic plane, as the set of points $(x,y)$ with $x\ge
0$. If we try to define the density of this region by circles all
expanding about a common center $c$, it is easy to see that the
density would depend on $c$, with any value strictly between 0 and 1
being obtainable for appropriate $c$.  This suggests that we will want
the origin, used for the expanding regions in (1), to be arbitrary.
\vs.1 \nd
{\bf Example 2 (stripe model).}
\vs.1
We now give a simple example of a region $S$ in the hyperbolic plane
such
that, when we try to define the density of $S$ relative to a sequence
of circles expanding about some point, we get the kind of oscillation
we found in the Euclidean annulus example.
We define the ``stripe model'' in the (upper half plane model of the)
hyperbolic plane, where the stripes are the regions separated by the
horocycles $h_j$, $j\in \Z$, defined by $y=y_j\equiv
e^{(j+1/2)\W}$, where fixed $\W> >1$ is to be specified. These curves
are equidistant by $\W$ in the hyperbolic metric. We call those
stripes separated by $h_{2j}$ and $h_{2j+1}$ ``black'', and the
others ``white'', and we declare the region $S$ of interest to be the
union of the black stripes.

Consider the circle with hyperbolic center $c=(0,1)$ and hyperbolic
radius $R=(N+1/2){\W}$, where $N>>1$ is to be specified. We will use the
following relations between the hyperbolic center $(H,K)$ and
hyperbolic radius $R$ of a given circle and its Euclidean center
$(h,k)$ and Euclidean radius $r$:

$$h=H,\ \ k^2-r^2=K^2,\ \ r=k\tanh(R).\eqno(8)$$

\nd So our circle has Euclidean center $(0,\cosh[R])$ and Euclidean
radius
$\sinh(R)$.

We will show that, if $N$ is even, the area inside the circle, of
the black stripes is larger than that of the white stripes; in
particular, each black stripe, between $h_j$ and $h_{j+1}$, $j\le
N-3$, is larger (by a factor 2) than that of the neighboring white
stripe above it (between $h_{j+1}$ and $h_{j+2})$, and therefore the
area of the circle is at least 2/3 black.

For $-N-1\le j\le N-1$, the area $A_j$ of the stripe between $h_j$ and
$h_{j+1}$ is:

$$\eqalign{
A_j&=\int_{y_j}^{y_{j+1}}\int_{-[2y\cosh(R)-1-y^2]^{{1\over 2}}}
^{[2y\cosh(R)-1-y^2]^{{1\over 2}}}{1\over y^2}\,dx\,dy\cr
&=\int_{y_j}^{y_{j+1}}{2[2y\cosh(R)-1-y^2]^{{1\over 2}}\over y^2}\,dy
.}\eqno(9)$$

\nd For $-N\le j \le N-2$ the leading behavior
as $N,\ \W\to \infty$ (and recalling that $R=[N+1/2]\W)$, is

$$\eqalign{
A_j&\sim \int_{y_j}^{y_{j+1}}{2y^{1/2}e^{R/2}\over y^2}\,dy \cr
&\sim 4e^{[R/2-(j+1/2)\W/2]}}\eqno(10)$$

\nd where $a\sim b$ means $\displaystyle {a\over b}\to 1$ as $N,\ \W\to
\infty$.
So

$${A_j\over A_{j+1}}\sim e^{\W/2}.\eqno(11)$$

\nd For $j=-N-1$ we have:

$$\eqalign{
A_{-N-1}=&\sim \int^{e^{-R+\W}}_{e^{-R}} {2(e^Ry-1)^{1/2}\over
y^2}\,dy\cr
&\sim 2e^R\int_1^{e^\W}{(z-1)^{1/2}\over z^2}\,dz\cr
&\gtrsim 2e^R\int_2^{e^\W}{1\over z^2}\,dz\cr
&\gtrsim e^R}\eqno(12)$$

\nd so

$${A_{-N-1}\over A_{-N}}\gtrsim {1\over 4}e^{\W/2}.\eqno(13)$$

\nd Finally we note that ${1\over 4}e^{\W/2}$ can be made as large as
desired,
in particular larger than 2, which completes the argument that the
relative densities of the set $S$ of black stripes does not have a
well defined limit.

\vs.1

The example of the stripe model in the hyperbolic plane, where the
stripes are all of equal ``width'', is more unsettling than the
example of annuli in Euclidean space discussed above, where in a sense
the oscillation was more obviously built in. We will see below that
this stripe model is only a simple version of a well known disk
packing.

There has been another common way to compute or estimate the density
of packings in Euclidean spaces, using tilings associated with the
packings, and the relative densities of the bodies in the tiles. (A
tile is a homeomorphic image of the closed unit ball, and a tiling is a
packing by tiles for which the union of the tiles is the full space
$\S$.) We emphasize that this is an attempt to reduce the intuitive
global idea of density, which involves taking a limit of approximate
densities in expanding regions of finite volume, to a more local
notion.  As a significant example of this approach we note an elegant
proof [Fe1, Rog] of the optimal density for packings of
equal disks in the Euclidean plane. The proof uses the Vorononi cells
of the bodies of a packing, where the cell for a body $\beta$ is
the set of all points $p\in \S$ as close to $\beta$ as to any body
of the packing. The proof shows that the
relative density in its Voronoi cell of any disk of any packing is
bounded above by that of any of the Voronoi cells in the obvious
hexagonal packing. This argument was extended to sphere packings in
$\S$ by K. B\"or\"oczky, who showed [Bo2] that the relative density of
any sphere of any packing of $\S$ in its Voronoi cell is bounded above
by the relative density associated with that of a regular simplex.
(See [FeK] for details.) Such relative
densities in tiles of associated tilings have remained an important
tool in analyzing optimal densities of sphere packings in Euclidean
spaces [FeK, Bez].

\vs.1 \nd
{\bf Example 3 (tight radius packings)}
\vs.1
In hyperbolic space, particularly the plane $\H^2$, the above method
of estimating or computing a density of sphere packings through an
associated tiling has been used convincingly for the special case of
disks of ``tight'' radius. The radius $r$ of a sphere in $\S$ is
called tight if the regular simplex of side length $2r$ admits a
(full-face to full-face) tiling of $\S$. In $\H^2$ this is the case if
and only if the equilateral triangle of edge length $2r$ has angles of
the form $2\pi/n$ for some $m\ge 7$, in which case

$$2r =2r_m=\cosh^{-1}[\cot({\pi\over m})\cot({2\pi\over m})]\eqno(14)$$

\nd and clearly $r_m\to \infty$ as $m\to \infty$. For disks with
tight radius $r_m$ the obvious ``periodic'' packing, in which each
disk is surrounded by $m$ disks touching it, has a well defined
density in the sense that, besides the method using Voronoi
tilings, {\it any} reasonable way to compute the density would give
the same value (namely $[3 \csc({\pi/m})-6]/[m-6]$ [Fe5]),
in particular any limit of the form (1) [BoR].

\vs.1 \nd
{\bf Example 4 (B\"or\"oczky's packing).}
\vs.1
There is an influential example due to B\"or\"oczky [Bo1] which
points out a difficulty in using relative density in tiles to define
the density of at least some packings in hyperbolic space, even some
which are rather symmetric. Place disks in the upper half plane model
of the hyperbolic plane with Euclidean centers at those points with
coordinates

$$ \{(e^{2j+{1\over 2}}(k+{1/2}),\ e^{2j+{1\over 2}})
\,|\,j,k\in\Z\}.\eqno(15)$$

\nd The connection between this and the stripe model is simple: we are
placing the disks equally spaced in the black stripes (and we are
taking the value $W=1$ for the width of the stripes). See Figure 1 for a
picture of the packing, which includes some horocycles and geodesics
to help understand the structure.

In Figure 2 we see the same packing with two congruent tiles in dark
outline. For each tile consider the tiling of the plane made by
congruent copies of the tile, as follows. First produce copies of the
tile by the congruences: $(x,y)\to (x + mw, y)$, $m\in
\Z$, where $w$ is the Euclidean width of the body. This fills out
a black and white stripe. Then produce, from these, more copies of the
tile by the congruences: $(x,y)\to (e^{2m}x,e^{2m}y)$, $m\in
\Z$. Together these copies of the original tile will cover the whole
plane. The two tilings made this way, one from each of the tiles in
Figure 2, are both simply related to the same packing of disks. The
punchline is, the tiling made by starting with the tile on the left in
Figure 2 would suggest assigning a ``density'' of the packing of disks
twice the value suggested by the tiling made by starting with the tile
on the right! We repeat the point that using a tiling to compute the
density of some packing, thus making the computation more local, is
useful in Euclidean spaces but is less convincing in hyperbolic
spaces.

\vs.1

We now return to the question of a definition of optimally dense
packings of $\H^n$. As we say above, for packings of Euclidean space
the notion of densest packings is easy to clarify, and one way to
understand this is through the computation of the ratio $f(\rho,a)$ of
volumes of concentric spheres of radii $\rho$ and $\rho+a$.

Note that:

\nd i) in $\E^n$ $\displaystyle f(\rho,a)\equiv {\rho^n\over
(\rho+a)^n}$,
so for fixed $a>0$ and $n$,
$f(\rho,a)\to 1$ as $\rho\to \infty$;

\nd ii) in $\E^n$, for fixed $a>0$ and $\rho$,
$f(\rho,a)\to 0$ as $n\to \infty$;

\nd iii) in $\H^n$, for fixed $a>0$ and $n$,
$f(\rho,a)\to e^{-ca}$ as $\rho\to \infty$, for some constant $c>0$.

To see why these phenomena interfere with a generalization to hyperbolic
space
of the method used earlier for Euclidean packings, consider the
packings $\P_\rho$ of the hyperbolic plane, by disks of fixed radius
$R$, defined for each $\rho>>0$ as follows. For each sufficiently
large radius $\rho >>R$, place disks of radius $R$ on the
circumference of a circle $C_\rho$ of radius $\rho$, so that: they
cover all but perhaps one arc of the circumference; there are as many
disks as possible without overlap; disks intersect only at points of
the circumference. We now show that by taking $R$ (and therefore
$\rho$) large enough we can ensure that the fraction of the area of
$C_\rho$ covered by the disks is as close to 1 as desired.

The fraction of the area of $C_\rho$ which is in the annulus between
$C_\rho$ and the concentric $C_{\rho'}$ for $\rho > \rho'$ is of the
order $1-e^{\rho'-\rho}$ for large $\rho,\ \rho'$, and by taking $0<<
\rho - \rho'<< R << \rho' << \rho$ we can ensure that most of this
area is inside the disks of radius $R$ -- all except those regions
outside pairs of touching disks of radius $R$ and outside the circle
$C_{\rho'}$, plus the region near any uncovered arc of
$C_\rho$.  But using the convexity of circles, the former regions are
each contained in triangles of the form $TUV$ (see Figure 7), so have
negligible
area, and another simple triangle argument applies to the region near
any uncovered arc of $C_\rho$.

So by choosing $R$ appropriately we could get almost all the area
of $C_\rho$ to lie outside $C_\rho'$.

Where in the Euclidean argument we used larger and larger cubes, in
hyperbolic space we would use fundamental domains of cocompact
subgroups of the isometry group $\G$ of $\H^n$. But we needed the fact,
in Euclidean space, that the volume of the portion of a
packing near the boundary of the fundamental domain would be
negligible, while we see now that for large fundamental domains
and large bodies, this is far from the case.
In summary, where we used i) to show the existence of optimal packings
in Euclidean space, in hyperbolic space we have instead iii), which
for large spheres is approximately ii). This is the intuitive reason
why there has been difficulty defining optimally dense packings
in hyperbolic space for so long.
\vs.2 \nd
{\bf II. Some responses to the problem}
\vs.1
We have summarized above the arguments that, for at least some
packings in hyperbolic space, there seems to be no reasonable notion
of density. We now consider how one might proceed with an analysis of
optimal density.

Two avenues of response that come to mind are: to replace the
essentially global definition of density with something more local; or
to find a way to define density for those packings where it is
reasonable, together with a convincing argument for excluding the
others.
\vs.1 \nd
{\bf a. Completely Saturated Packings}
\vs.1

A packing is called completely saturated [FeK, FKK, Kup]
if it is not possible to
replace a finite number of bodies of the packing with a greater total
volume
of
bodies and still remain a packing.  Intuitively, we think of a
completely
saturated packing as one that is locally densest. In [FKK], it was
proven
that any convex body of Euclidean space admits a
completely saturated packing (and more generally any body with the
strict
nested similarity property) (see also [Kup]). In [Bow], it is proven
that
completely saturated packings exist for all bodies $\beta$ in either
Euclidean of hyperbolic space. The argument given there extends easily to
finite collections $\B$.
\vs.1 \nd
{\bf Example 5 (a completely saturated packing with low density).}
\vs.1
As pointed out in [FKK], a completely saturated packing of Euclidean
space
is a densest packing. This is not true in hyperbolic space. In this example,
we construct a pair of bodies $\beta_1,\beta_2$ in
$\H^2$ such that two completely saturated periodic packings by
$\{\beta_1,
\beta_2\}$ exist that have different densities. The
reason, as we will see, is due to the fact that the length of the
boundary of a region in the hyperbolic plane is comparable to its area.
Let $\beta_1$ be the tile shown in Figure 8. It is a regular octagon
with all interior angles equal to $2\pi/8$. Let $\T_1$ be the unique
periodic tiling by $\beta_1$. Let $\beta'_2$ be the tile shown in Figure 9.
It is formed from $\beta_1$ by adding ``protrusions'' to some edges
and
``indentations'' to others. We will assume that these protrusions and
indentations are made so that they fit together but are narrow enough so
that there is a region of finite area $C_1$ in each indentation that
cannot be occupied by a nonoverlapping copy of $\beta'_2$ unless it is
occupied by a protrusion. Also we assume that each protrusion fits into
a
unique indentation.

$\beta'_2$ admits a unique periodic tiling $\T_2$. Let $\beta_2$ be
equal
to $\beta'_2$ with a small hole removed from its interior. Let $\P$ be
the
obvious periodic packing by $\beta_2$ (i.e. the one that comes from
$\T_2$
by removing a small hole from the interior of each tile). Since $\beta_1$
admits a periodic tiling, it is clear that the optimal
density of $\{\beta_1,\beta_2\}$ is one. Just as clear, is the fact that
the density of $\P$ is $area(\beta_2)/area(\beta_1) < 1$. We will show
that $\P$ is completely saturated (if the hole in $\beta_2$ is small
enough).

It is a standard fact of hyperbolic geometry that there exists a
constant
$C_2>0$ depending only on the symmetry group of $\T_1$ (and the fact
that
$\beta_1$ contains a fundamental domain for this group) such that for
all
finite subtilings $\T'$ of $\beta_1$, $|\partial \T'| \ge C_2|\T'|$ (by
$|\partial \T'|$ we mean the number of edges contained in exactly one
tile
of $\T'$ and by $|\T'|$ we mean the number of tiles in $\T'$). Since the
hole in the interior of $\beta'_2$ can be made as small as we like, we
may
assume that $area(\beta_2) > area(\beta_1) - C_1C_2/2$.

Suppose for a contradiction that $\P$ is not completely saturated. Then
there exists a finite subpacking $\P' \subset \P$ and another finite
packing $\P''$ such that $(\P - \P') \cup \P''$ is a packing and
$area(\P'') > area (\P')$. We may assume without loss of generality that
$\P \cap \P'' = \emptyset$.

We claim that the number of edges of $\P'$ that have protrusions on them
coming from bodies of $\P'$ is at least $|\partial \P'|/2$. So let $e$ be
any edge on the boundary of $\P'$. Let
$e=e_0,e_1,..,e_n$ be the sequence of edges defined by for $1 \le i <
n$,
$e_{i+1}$ and $e_i$ are on a body of $\P'$ and $e_{i+1}$ is the
``opposite
side'' of $e_i$ in the sense that if $e_i$ has a protrusion on it
(relative
to the body containing both $e_i$ and $e_{i+1}$) then $e_{i+1}$ is its
corresponding indentation and vice versa. This sequence is uniquely
defined and ends in an edge $e_n$ on the boundary of $\P'$. It is easy
to
see that if $e_0$ corresponds to an indentation of $\P'$ (i.e. $e_0$ has an
indentation on it coming from a body of $\P'$) then $e_n$
corresponds to a protrusion and vice versa. Thus the claim is proven.

Note that it is not possible for any body of $\P''$ to fill completely
any indentation on the boundary of $\P - \P'$ (in fact a region of area
at
least $C_1$ is always unfilled). Hence the total area of $\P''$ is at
most

$$\eqalign{
area(\P'') \le& |\P'| area(\beta_1)  - (C_1/2)|\partial \P'| \cr
               \le& |\P'|[area(\beta_1) - C_1C_2/2]\cr
               <& |\P'|area(\beta_2)\cr
               =& area(\P').}\eqno(16)$$

This contradicts the choice of $\P''$. So $\P$ is completely
saturated. The moral is that, in hyperbolic space, locally densest does
not imply globally densest.

\vs.1 \nd
{\bf b. Controlling pathological packings}
\vs.1
We now discuss an approach to density specifically aimed at
controlling those packings, such as the above example of B\"or\"oczky,
which pose difficulty in computing a reliable density. Even though the
methods are also applicable to Euclidean space, the interests of this
article make it natural to specialize the discussion from now on to
$\S=\H^n$.

The key idea is to use a pointwise ergodic theorem of Nevo ([Nev,
Thm. 1] for dimension $n\ge 3$; [NeS, Thm. 3] for $n\ge 2$), the
conclusion of which is the existence of limits of the type (1) in the
intuitive definition of density. The fact that such theorems only
prove existence of the limit ``almost everywhere'' is not a defect, it
is a feature, necessitated by examples such as that of B\"or\"oczky.

We begin by reproducing some notation and results from [BoR].
Let $d(\cdot,\cdot)$ be the usual metric on $\S$, and let $\O$ be a
distinguished origin. We suppose given a finite collection $\B$
of bodies $\beta_j$ in $\S$. Let $\Sigma_\B$ be the space of all
``relatively-dense'' packings of $\S$ by congruent copies of the $\beta_j$,
that is, packings $\P$ with the property that any congruent copy of a
body in $\B$ intersects a body of $\P$. On $\Sigma_\B$ we put the
following metric, corresponding to uniform convergence on compact
subsets of $\S$:

$$d_\B(\P_1, \P_2)=\sup_{k\ge 1}{1\over k}h(B_{k}\cap \P_1,B_{k}\cap
\P_2),\eqno(17)$$

\nd where ${B}_{k}$ denotes the closed ball of radius $k$ centered
at the origin, and for compact sets $A$ and $C$ we
use the Hausdorff metric

$$h(A,C) \equiv \max \{ \sup_{a \in A}\inf_{c\in C}d(a,c), \sup_{c
\in C}\inf_{a\in A}d(a,c)\}.\eqno(18)$$

\nd It is not hard to see [RaW] that $\Sigma_\B$ is compact in this
metric
topology, and that the natural action: $(g,\P)\in \G\times
\Sigma_\B\longrightarrow g(\P)\in \Sigma_\B$ of the isometry group $\G$
of
$\S$ on $\Sigma_\B$ is (jointly) continuous.
Let $\M(\B)$ be the family of Borel
probability measures on $\Sigma_\B$. We call a measure $\mu \in \M(\B)$
``invariant'' if for any Borel subset $E \subset \Sigma_\B$ and any $g \in
\G$,
$\mu(gE) = \mu(E)$. Let $\M_I(\B)$ be the subset of invariant measures
and $\M^e_I(\B)$ the convex extreme (``ergodic'') points of $\M_I(\B)$,
all in their weak* topology, in which $\M(\B)$ and $\M_I(\B)$ are
compact.

We will study these ergodic measures as a substitute for studying individual
packings. As we will see, for any ergodic measure $\mu \in \M_I(\B)$ there
is a set of packings $Z$ of full $\mu$-measure such that for each $\P \in
Z$, the orbit of $\P$ is dense in the support of $\mu$. So studying $\mu$ is
a lot like studying a packing in  $Z$. We will make this relationship more
clear in what follows but first some examples.

Suppose $\P$ is a ``periodic'' packing, i.e.\ the symmetry
group
$\Gamma_\P$ of $\P$ is cocompact in $\G$. We will construct a measure
$\mu_\P\in \M^e_I(\B)$ whose support is contained in the orbit
$O(\P)\equiv\{g\P\, | \, g\in \G\}\subset \Sigma_\B$ of $\P$.  $O(\P)$
is
naturally homeomorphic to the (metrizable) space $\G/\Gamma_\P$ of left
cosets by the homeomorphism $q_\P: O(\P)\to \G/\Gamma_\P$ with
$q_\P(g\P)=g\Gamma_\P$.  There is a natural probability measure on
$\G/\Gamma_\P$ induced by Haar measure on $\G$ by the projection map
$\pi_\P: \G \to \G/\Gamma_\P$.  (Aside from an overall normalization the
measure on $\G/\Gamma_\P$ can be defined on sufficiently small open
balls $B\subset \G/\Gamma_\P$ as the Haar measure of any of the
components of $\pi_\P^{-1}(B)$.)  Hence $q_\P$ induces a probability
measure $\hat \mu_\P$ on $O(\P)$. This measure can then be extended to
all of $\Sigma_\B$ in the following way: $\mu_\P(E) = \hat \mu_\P[E \cap
O(\P)]$ for any Borel set $E \subseteq \Sigma_\B$. We will use the term
``periodic measure'' to denote any measure in $M_I(\B)$ associated in
this way with the orbit of a periodic packing.  It is not hard to
prove from the uniqueness of Haar measure on $\G$ that there is only
one probability measure, with support in the orbit of a periodic
packing, which is invariant under $\G$.

Next, we define the density of an invariant measure. After the definition,
we will show how the density of an invariant measure relates to the density
of packings in its support.

For any $p\in \H^n$ we define the real valued function $F_p$ on
$\Sigma_\B$ as the indicator function of the set of all packings $\P$
such
that $p$ is contained in a body of $\P$. (The latter condition will
sometimes be expressed as $p\in \P$.)
\vs.1 \nd {\bf Definition 1}. For any invariant measure $\mu\in
\M_I(\B)$, the ``average density'' $D(\mu)$ is defined as
$\int_{\Sigma_\B}
F_p(y)\, d\mu(y)$.
\vs.1 \nd Note: the average density $D(\mu)$ is
independent of the choice of $p$, because of the invariance of the
measure, so $p$ is not needed in the notation. For convenience we
sometimes use $p=\O$.

If $\P_\mu$ is a random packing with distribution $\mu$ then the above
definition states that the density of $\mu$ is the probability that the
origin is contained in a body of $\P_\mu$.

For periodic packings $\P$ there is an obvious notion of density using
a fundamental domain of $\Gamma_\P$. The above definition
of density coincides with this intuitive notion for such special $\P$.
\vs.1 \nd {\bf Proposition 1 [BoR]}. If $\P$ is a periodic packing,
$D(\mu_x)$ is the relative volume of any fundamental domain for
$\Gamma_x$ taken up by the bodies of $x$.
\vs.1

We need the following notation. As usual we let $\G$ denote the group of
orientation preserving
isometries of
hyperbolic $n$-space $\H^n$ (for some fixed $n \ge 2$). Let
$\pi: \G \to \H^n$ be the projection map $g \to g\O$ where $\O$ is
some distinguished point in $\H^n$. Then we let $\tilde B_r$ denote the
inverse image under $\pi$ of the closed ball of radius $r$ centered at
$\O$. Finally let $\lambda_\G$ denote a Haar measure on $\G$, normalized
so
that $\lambda_\G(\tilde B_r)$ is the volume of the $r$-ball in
$\H^n$.

We will use the following special case of Theorem 3 in [NeS] to relate the
density of an ergodic measure to the density of (almost every) packing in
its support.

\vs.1 \nd
{\bf Theorem 1 [Nevo]}. If $\G$ acts continuously on a compact metric
space $X$ such
that there is a Borel probability measure $\mu$ on $X$ that is invariant
and ergodic under this action, then for every function $f \in
L^p(X,\mu)$
($1 < p < \infty$) there is a set $Z$ of full $\mu$ measure such that
for
every $z \in Z$,

$$
\int_X f d\mu = \lim_{r \to \infty} {1\over \lambda_\G(\tilde B_r)}
\int_{\tilde B_r} f(gz) \, d\lambda_\G(g).\eqno(19)
$$

\noindent Actually we will use the following extension of this result.

\vs.1 \nd
{\bf Theorem 2}. Under the same hypotheses as the above theorem, the set
$Z$
may be taken to be invariant under $\G$.
\vs.1

We will prove this result in the next section. Applying Theorem 2 to the
function $F_p$ and using the proof of Prop.\ 2 of [BoR] we get

\vs.1 \nd
{\bf Theorem 3}. If $\mu \in M^e_I(\B)$ then there exists a set
of packings $Z$, of full
$\mu$-measure, such that for all $p \in \H^n$ and all $\P \in Z$

$$\lim_{r\to \infty} {vol[\P\cap B_p(r)]\over vol[B_p(r)]} =
D(\mu).\eqno(20)$$
\vs.1

Note that this implies that the (closure of the orbit of the) stripe model
has measure zero
with respect to every invariant measure. We will give another explanation
for this fact in a later section.

  From example 4 we concluded that it is not possible, in general, to
compute the density of a hyperbolic packing using a certain tiling 
associated to it.
In spite of this we will show that it is possible to compute the density of
an invariant measure using an associated space of tilings.

Let $\Sigma$ be a (compact, invariant) space of packings of $\H^n$. Let
$\mu$ be a $Isom^+(\H^n)$ invariant measure on $\Sigma$. Suppose that $\t$
is a space of tilings (of $\H^n$) and that there is an equivariant map
$\phi: \Sigma \to \t$. For example, $\t$ may be the space of Voronoi
tilings [FeK]
corresponding to $\Sigma$. For $\P \in \Sigma$ such that the
origin is contained in a tile of $\phi(\P)$, let $\tau_p(\P)$ denote the
tile of $\phi(\P)$ containing the point $p \in \H^n$. We claim that for any
$p$:

$$D(\mu) = \int_\Sigma \, {vol[\P \cap \tau_p(\P)] \over vol[\tau_p(\P)]} \,
d\mu(\P).\eqno(21)$$

Define a function $f: \H^n \times \H^n \times \Sigma \to \R$ by $f(p,q,\P) =
1/vol[\tau_p(\P)]$ if $p$ and $q$ are both in the tile $\tau_p(\P)$ and
$p$ is in a body of $\P$ (otherwise $f(p,q, \P) = 0$).
Define a measure $\nu$ on $\H^n \times \H^n$ by

$$\nu(E \times F) = \int_{\Sigma} \int_{E} \int_{F} f(p,q,\P) \, dvol(p) \,
dvol(q) \, d\mu(\P).\eqno(22)$$

Since $\mu$ is invariant, it is easy to check that for all $g \in
Isom^+(\H^n)$, $\nu(gE \times gF) = \nu(E \times F)$. The mass-transport
principle [BeS] implies that $\nu(E \times \H^n) = \nu (\H^n \times E)$ for
any measureable $E \subset \H^n$. But it can easily be checked that $\nu(E
\times \H^n) = vol(E)D(\mu)$ and $\displaystyle \nu(\H^n \times E) =
vol(E)\int_\Sigma \,
{vol[\P \cap \tau_p(\P)]/ vol[\tau_p(\P)]} \, d\mu(\P)$ (for any $p$).
This proves the claim. For emphasis, we repeat that when $\mu$ is an
invariant measure we can compute its density with respect to local
structures such as the Voronoi tilings. If $\mu_\P$ is a periodic
measure and $\t$ is the space of tilings by a fundamental domain of
$\Gamma_\P$ then there is a natural equivariant map from the orbit of $\P$
to $\t$. The above result then yields Proposition 1.

We now define optimality through measures.
\vs.1 \nd {\bf Definition 2}. $D(\B)\equiv \sup_{\mu\in
\M^e_I(\B)}D(\mu)$
will be called the ``optimal density for $\B$'', and any ergodic
measure $\tilde \mu\in \M^e_I(\B)$ will be called ``optimally dense
(for $\B$)'' if $D(\tilde \mu)=D(\B)$. We define ``optimally dense
packings'' a little differently than in [BoR]. We say that a packing $\P$ is
optimally dense if there is an optimally dense measure $\mu$ such that the
orbit of $\P$ is dense in the support of $\mu$ and for every $p \in \H^n$
$D(\mu)$ is equal to the limit of the relative fraction of volume in
expanding spheres centered at $p$ taken up by bodies of $\P$.

\vs.1
One of the main results of [BoR] asserts the existence of
optimally dense measures. Also, it was proven that for every $\mu \in
M^e_I(\B)$ there exists a set $Z$ of  full $\mu$-measure such that for every
$\P \in Z$, the orbit of $\P$ is dense in the support of $\mu$. Using
Theorem 3 this implies

\vs.1 \nd {\bf Theorem 4}. For any finite collection $\B$ of bodies
there exists an optimally dense measure $\mu$ on $\Sigma_\B$, and a
subset of the support of $\mu$, of full $\mu$-measure, of optimally
dense packings.
\vs.1 \nd
Note: There may be many optimally dense measures for a given $\B$.
\vs.1

In [Bow] it is proven that the set of completely saturated packings has
full-measure with respect to any optimally dense measure $\mu$. In other
words, if $\P_\mu$ is a random packing with optimally dense distribution
$\mu$ then $\P_\mu$ is completely saturated almost surely. Example 5 shows
that the converse is false.

A major result of [BoR] was that the set of all radii $r$ such that there
exists an optimally  dense {\it periodic} measure for the sphere of radius
$r$ (in hyperbolic space) is at most countable. Thus most
optimally dense sphere packings are complicated.

\vs.1 \nd
{\bf c. What Invariant Measures Avoid}
\vs.1

Some packings, such as the B\"or\"oczky example, do not have a well-defined
density. We claim this is ``due to'' the fact that the closure of the orbit
of
such a
packing has measure zero with respect to every invariant measure $\mu$. In
this section we prove this statement and show other examples of packings
that are not ``seen'' by invariant measures.

\vs.1 \nd
{\bf Example 6 (Penrose's binary tilings)}
\vs.1
Perhaps the most relevant to the discussion in section I is the $\B$
consisting of the body $\beta$ shown in Figure 3. (This is a minor
variation on the tile in [Pen], and a special case of tiles in [MaM].)
We know copies of this
body can tile $\H^2$, and since limits in $\Sigma_\B$ of tilings will
again be tilings, if there were any invariant measure $\mu\in
\M_I(K)$ with support in the orbit closure of such a tiling it would
clearly have density 1. However we can see there is no such measure as
follows. First consider the slightly simpler, and better known,
example of the natural action of the isometry group $\G$ of $\H^2$
(namely $\G=PSL_2(\R)$) on the boundary $\Delta$ of $\H^2$, instead
of its action on the set of tilings. Assume
there is a measure $\mu$ on $\Delta$ invariant under $\G$. Any
hyperbolic element $g_h\in \G$ has 2 fixed points in $\Delta$, \
$p_1,\ p_2$, and moves all other points towards one and away from the
other. From its invariance under $g_h$, $\mu(\{p_1,p_2\})=1$.  Then
considering that any elliptic element $g_e\in \G$ has no fixed points
in $\Delta$, and $\mu$ must also be invariant under $g_e$, we get a
contradiction. So there are no probability measures on $\Delta$
invariant under $\G$. Going back to our space $\t_\B$ of tilings by
our body $\beta$, consider the function $f$ from $\t_\B$ to $\Delta$,
which takes each tiling to the point ``pointed to'' by the protrusion
on each body in the tiling. $f$ is obviously continuous. If there were
a probability measure $\mu$ on the space $\t_\B$ of tilings,
invariant under the action of $\G$, we could define a corresponding
measure $\mu_f$ on $\Delta$ by $\mu_f(E)= \mu(f^{-1}[E])$. Since no
such $\mu_f$ exists, this proves no such $\mu$ exists.

Now assume the optimal density for $\B$, $D(\B)$, is 1, with
an optimal measure $\mu$. For each $R>0$ consider the function
on $\Sigma_\B$

$$f_R(\P)={1\over vol(B_R)}\int_{{\tilde B}_R}F_\O(g\P)\,
d\lambda_\G(g),\eqno(23)$$

\nd which gives the relative area of the ball $B_R$ covered by the disks
of $\P$. From the invariance of of $\mu$,

$$\int_{\Sigma_\B}f_R(\P)\,d\mu(\P)=
\int_{\Sigma_\B}F_\O(\P)\,d\mu(\P)=D(\mu)=1,\eqno(24)$$

\nd so $f_R(\P)=1$ for $\mu$-almost every $\P\in \Sigma_\B$. Letting $R$
run through the positive integers, and intersecting the sets of full
measure we get for each such $R$, we see there is a set of packings of
full measure which are tilings. Since the closure of a set of tilings
can only contain tilings, and the support of $\mu$ must be invariant
under $\G$, that support is contained in the set of all tilings of $\beta$.
But we saw above that there can be no such measure as $\mu$, and this
proves that $D(\B)\ne 1$.
(Using modifications of this
example
it can be shown [Bow] that for every $\epsilon > 0$ there exists a body
$\beta$ that admits a tiling of $\H^n$ and $D(\beta) < \epsilon$.)

\vs.1

The formalism above leads one to assert that the densest
packings of the body $\beta$ have density bounded away from 1, even
though
one can tile $\H^2$ with copies of $\beta$. The ``reason'' for this is
that there are no invariant measures which can ``see'' the tilings;
they are a set of measure zero for every invariant measure on the
space of all packings by $\beta$. We explore the consequences of
this using some of the examples we discussed earlier.

Consider again the tile $\beta$ shown in Figure 3. Congruent copies of
$\beta$
can only tile the plane (up to an overall rigid motion) as
in Figure 4 (in which the little bumps on the tiles are not shown.)
Construct the tile $\bar \beta$ of Figure 5 out of three abutting copies
of $\beta$.  Now drill a hole in $\bar \beta$, producing the body
$\bar \beta_0$, as shown in Figure 6. Note that the packings of the
plane by $\bar \beta_0$ obtained in the obvious way from the tilings
by $\bar \beta$, are precisely the complements of the disk packings of
B\"or\"oczky discussed above. The point is, although it might seem
reasonable to assign a density of 1 to the tiling of Figure 4, that
would seem to imply a well defined density to the packing of Figure 2,
which we know is misleading. In other words, the meaningfulness of the
density of the tiling of Figure 4 is unstable under arbitrarily small
perturbations (drilling arbitrarily small holes). Notice that when we
drill these small holes we turn the tiling into a mere packing,
forcing us to give up the ``simplicity'' of the tiling, as a global
object with seemingly obvious density, and leaving us to find some
meaningful way to assign a density to the resulting packing.  As we
will see below, the difficulty in assigning a density to a packing,
for instance congruent copies of a single body $\beta$, can derive
from the complexity of the set of rigid motions of $\beta$ that define
the packing. And in this sense a tiling is no simpler; treating it as
a global object with an ``obvious'' density simply avoids coming to
grips with the essential nature of the assignment of density for
packings.

In other words, the phenomenon whereby the ``optimal'' density can be
less (even far less) than 1 for a body which can tile space, can be
understood as related to the instability of the meaningfulness of the
density of the tilings under removal of small holes in the tiles. This
suggests that even for tilings one needs to keep track of the
individuality of the tiles. In this example that amounts to noting the
various sets of congruences used in producing the tilings; in some
sense those sets of congruences are too complicated to be analyzed
through our density formalism.

We have shown that if $\T$ is a tiling whose orbit closure in the space
of
packings factors onto the space at infinity of the hyperbolic plane,
then
there are no invariant measures on the orbit closure of $\T$. All of our
examples of strange behaviour in the hyperbolic plane have, so far, been
constructed using this principle. Could this be the only way of
constructing such examples?
\vs.1 \nd
{\bf Example 7.}

The following example is a variant on a simple construction. First consider
two congruent regular all-right-angles octagons in the plane. Label their
edges in clockwise order $e_1, e_2,..,e_8$ and $e'_1, e'_2,..,e'_8$. By
identifying $e_k$ with $e'_k$ for odd $k$ (by orientation reversing
homeomorphisms), we obtain a sphere $X'$ with four open disks removed. If we
then identify pairs of boundary components of $X'$, the resulting object is a genus two surface. The covering
space $S$ of the genus two surface corresponding to the commutator subgroup
looks like the boundary of a regular neighborhood of the standard Cayley
graph of the free group on two generators. In this way, we obtain a tiling
of $S$ by regular all-right octagons.

For the variation, we wish to distinguish a boundary component of $X'$. We
do this by modifying the edges of the all-right octagon so that $e_2$ (and
$e'_2$) has a protrusion and $e_4, e_6, e_8$ (and $e'_4, e'_6, e'_8$) have
indentations. We call this tile $\tau$ (see figure 10). We identify $e_k$
with $e'_k$ for $k$ odd as before to obtain $X$ as in figure 11. We want to
obtain a tiling of $S$ by $X$. For this picture the standard Cayley graph of
the free group on two generators. Draw arrows on each of the edges so that
each vertex has exactly one outgoing arrow. For each vertex $v$, let $X_v$
be a copy of $X$. If there is an edge with arrow pointing from $v$ to $w$,
then we identify the boundary component of $X_v$ that has the protrusion
with one of the boundary components of $X_w$ that has an indentation. In
this way, we obtain a tiling of $S$ by $X$ (and thus by $\tau$).

Note that the free group $F_2$ on 2 generators $\{a, b\}$ acts naturally
and isometrically on $S$. This action lifts to an isometric action of
$\H^2$ via the covering map. Though not relevant to what follows, note
that this lift is unique up to postcomposition by rigid motions of $\H^2$.

If we start a walk in $\T_S$ from some initial tile and follow the
protrusions we get ``closer'' to a point on the ideal boundary  of $S$.
It
is not too hard to see that this point does not depend on the initial
tile
chosen but only on the tiling $\T_S$. Therefore, there is a map from the
space of tilings of $S$ by $X$ (defined similar to the same way
$\Sigma_\B$ is defined)  to the ideal boundary of $S$ that commutes with
the action of $F_2$. Since the ideal boundary does not admit an
invariant
Borel probability measure (for practically the same reason that $\Delta$
does not admit an invariant measure), neither does the space of tilings
of
$S$ by $X$.

Suppose that there exists an invariant measure $\mu$ whose support is
contained in the orbit closure of $\T$. Then this measure pushes forward
via the covering map to a measure $\mu_S$ on the space of tilings on $S$
by $\tau$. This measure $\mu_S$ is invariant under the action of $F_2$
but
this contradicts the previous paragraph.

Now suppose that there is an equivariant map $\phi$ from the orbit
closure
$\overline{O(\T)}$ of $\T$ in $\Sigma_\tau$ to $\Delta$ the boundary at infinity of the hyperbolic plane. Let $p = \phi(\T)$.
Since
$\phi$ is equivariant, the stabilizer of $\T$ must be contained in the
stabilizer of $p$. However, the stabilizer of $\T$ is noncyclic (since
it
contains an isomorphic copy of the fundamental group of $S$ which is
noncyclic). By the theory of fuchsian groups, the stabilizer of $\T$ does
not fix any
point
at infinity. This contradiction shows that $\phi$ cannot exist.
\vs.2 \nd
{\bf III. Proof of Theorem 2}
\vs.1
We need the following well-known
fact (see 7.1.1.2 in [AVS]):

\vs.1 \nd
{\bf Lemma 2}. There exist positive real constants $c_1$ and $c_2$
(depending only on the dimension $n$)
such that
$$
\lim_{R \to \infty} vol[B_R]e^{-c_1 R} = c_2.\eqno(25)
$$

\vs.1 \nd
{\bf Corollary 1}.
For $r > 0$,
$$
\lim_{R \to \infty} {vol[B_{R-r}] \over vol[B_R]} = e^{-c_1 r}.\eqno(26)
$$

\vs.1 \nd

\nd If $h: \G \to \R$ is Borel, let $A_-(h), A_+(h): \G \to \R \cup \{\pm
\infty\}$ be
defined by
$$
A_-(h)(g) =\liminf_{R \to \infty} {1\over vol(B_R)} \int_{\tilde B_R}
h(g'g) \, d\lambda_\G(g')\eqno(27)
$$
$$
A_+(h)(g) =\limsup_{R \to \infty} {1\over vol(B_R)} \int_{\tilde B_R}
h(g'g) \, d\lambda_\G(g').\eqno(28)
$$

\vs.1 \nd {\bf Lemma 3}.
If $h$ is any nonnegative Borel function on $\G$ then $A_+(h)$ and $A_-(h)$
are continuous.

\vs.1 \nd Proof.

Let $g_1, g_2 \in \G$ be such that the distance between $g_1\O$ and $g_2\O$
is
$r$. Then

$$\eqalign{
A_+(h)(g_1) =& \limsup_{R \to \infty} {1\over vol(B_R)} \int_{\tilde B_R}
h(g'g_1) \, d\lambda_\G(g')\cr
            =& \limsup_{R \to \infty} {1\over vol(B_R)} \int_{\tilde B_R
g_1^{-1}} h(g') \, d\lambda_\G(g')\cr
            =& \limsup_{R \to \infty} {1\over vol(B_R)}
\int_{\pi^{-1}[B_R(g_1\O)]} h(g') \, d\lambda_\G(g')\cr
            \ge& \limsup_{R \to \infty} {1\over vol(B_R)}
\int_{\pi^{-1}[B_{R-r}(g_2\O)]} h(g') \,
d\lambda_\G(g') \cr
            =& \limsup_{R \to \infty} {vol(B_{R-r})\over vol(B_R)} {1 \over
vol(B_{R-r}) }
\int_{\pi^{-1}[B_{R-r}(g_2\O)]} h(g') \,
d\lambda_\G(g') \cr
            =& e^{-c_1 r}A_+(h)(g_2). \cr
}\eqno(29)$$

\noindent Since $g_1$ and $g_2$ are arbitrary, $A_+(h)$ is continuous. The
proof for $A_-(h)$ is similar.\qed
\vs.1 \nd
Proof of Theorem 2.
Let $\G_0$ be a countable dense subset of
$\G$. Let $Z_0$ be as in Theorem 1. At first we assume that $f$ is
nonnegative. Let $Z_f = \bigcap_{g \in \G_0}
g^{-1}Z_0$. Since $Z_f$ is a countable intersection of sets of measure 1,
$\mu(Z_f) = 1$. By definition, for all $z \in Z_f$ and for all $g \in \G_0$,
$gz \in Z_0$. For $x \in X$, define $h_x:\G \to \R$ by $h_x(g) =
f(gx)$. For $z \in Z_f$ and $g \in
\G_0$,
we have
$$\eqalign{
A_+(h_z)(g) =& \limsup_{R \to \infty} {1\over vol(B_R)} \int_{\tilde
B_R} h_z(g'g) \, d\lambda_G(g')\cr
=& \limsup_{R \to \infty} {1\over vol(B_R)} \int_{\tilde B_R}
f(g'gz) \, d\lambda_G(g')\cr
=& \int_X f d\mu.\cr
}\eqno(30)$$

The last equation holds since $gz \in Z_0$. By the previous lemma,
$A_+(h_z)$ is
continuous. So the above equations hold for all $g \in
\G$. Similarly, $A_-(h_z)(g) = \int_X f d\mu$ for all $g \in
\G$
and $z \in Z_f$. So

$$\lim_{R \to \infty} {1 \over vol(B_R)} \int_{\tilde B_R} f(g'gz) \, 
d\lambda_G(g') = \int_X f \, d\mu \eqno(31)$$

\nd for all $g \in \G$ and all $z \in Z_f$. So the set $Z = \cup_{g \in \G} 
gZ_f$ satisfies the conclusion of the theorem and completes the case when 
$f$ is nonnegative.

In general, we set $f=f_+ - f_-$ where $f_+$ and $f_-$ are nonnegative. By
the above, there are invariant sets $Z_+$ and $Z_-$ of full $\mu$ measure
satisfying the conclusion of the theorem for $f_+$ and $f_-$. The set $Z_f =
Z_+ \cap Z_-$ is invariant, of full $\mu$ measure and for all $z \in Z_f$,
we have
$$\eqalign{
\int_X f =& \int_X f_+ - f_- \, d\mu\cr
         =& \int_X f_+ \, d\mu - \int_X f_- \, d\mu\cr
         =& \lim_{R \to \infty} {1 \over vol(B_R)} \int_{\tilde B_R} 
f_+(gz)\, \lambda_G(g)\cr
         -& \lim_{R \to \infty} {1 \over vol(B_R)} \int_{\tilde B_R} 
f_-(gz)\, \lambda_G(g)\cr
         =& \lim_{R \to \infty} {1 \over vol(B_R)} \int_{\tilde B_R} f_+(gz)
- f_-(gz)  \, \lambda_G(g)\cr
         =& \lim_{R \to \infty} {1 \over vol(B_R)} \int_{\tilde B_R} f(gz)
\, \lambda_G(g).\cr
}\eqno(32)$$

\nd This proves the theorem.\qed
\vs.2 \nd
{\bf IV. Summary}

In [BoR] the notion of ``optimal density'' was
defined for packings in hyperbolic space $\H^n$, as above,
through the use of
probability measures, on a space of packings, invariant
under the congruence group of $\H^n$.
The notion of an optimally dense packing was also introduced in
[BoR], but not very successfully. The justification for
that term was not well connected to a limit (1); we could
only show there that for a set of packings of full measure, the
limit (1) existed relative to expanding spheres
centered about any countable set of centers.
One advance
in this paper is an extension of Nevo's ergodic theorem,
allowing us to extend this proof of existence to {\it all}
centers in $\H^n$, allowing a more natural notion of optimally
dense packings.

Perhaps more significantly, we have also contrasted our approach to
density with earlier approaches, and compared some key examples,
showing the significance of certain structural features of the space
of packings to the existence of well defined densities.

\vs.2 \nd
{\bf Acknowledgments}. We are grateful to Russell Lyons for pointing
out to us some uses of the Mass Transport principle and for finding several errors in a previous version of theorem 2..

\vfill \eject
\centerline{{\bf References}}
\vs.2 \nd
\item{[AVS]} D.V. Alekseevskij, E.B. Vinberg and A.S. Solodovnikov,
Geometry of spaces of constant curvature, in {\it Geometry II:
Spaces of Constant Curvature}, ed. E.B. Vinberg, Springer-Verlag,
Berlin, 1993.

\item{[BeS]} I. Benjamini and O. Schramm, Percolation in the
hyperbolic plane, {\it Jour. Amer. Math. Soc.} 14(2001) 487-507.

\item{[Bez]} K. Bezdek, Improving Rogers' upper bound for the density of 
unit ball packings via estimating the surface area of Voronoi cells from 
below in Euclidean d-space for all $d\ge 8$, {\it Discrete Comput. Geom.} 
28(2002), 75-106.

\item{[Bo1]} K. B\"or\"oczky, Gombkitoltes allando gorbuletu terekben, {\it
Mat. Lapok.} 25(1974), 265-306.

\item{[Bo2]} K. B\"or\"oczky, Packing of spheres in spaces of constant
curvature, {\it Acta Math. Acad.  Sci. Hung.} 32(1978) 243-261.

\item{[BoF]} K. B\"or\"oczky and A. Florian, Uber die dichteste Kugelpackung
in hyperbolischen Raum, Acta Math. Acad. Sci. Hung. 15(1964) 237-245.

\item{[Bow]} L. Bowen, On the existence of completely saturated packings and
completely reduced coverings, {\it Geometria Dedicata} (to appear).

\item{[BoR]} L. Bowen and C. Radin, Densest packing of equal spheres in
hyperbolic space, {\it Discrete Comput. Geom.} (to appear).

\item{[Fe1]} L. Fejes T\' oth, Uber einen geometrischen Satz, {\it Math. Z.}
46(1940) 79-83.

\item{[Fe2]} L. Fejes T\' oth, On close-packings of spheres in spaces of
constant curvature, Publ. Math. Debrecen 3(1953) 158-167.

\item{[Fe3]} L. Fejes T\' oth, Kreisaufuellungen der hyperbolischen Ebene,
Acta Math. Acad. Sci. Hung. 4(1953), 103-110.

\item{[Fe4]} L. Fejes T\' oth, Kreisuberdeckungender hyperbolischen Ebene,
Acta Math. Acad. Sci. Hung. 4(1953), 111-114.

\item{[Fe5]} L. Fejes T\'oth, {\it Regular Figures}, Macmillan, New York,
1964.

\item{[FeK]} G. Fejes T\'oth and W. Kuperberg, Packing and covering
with convex sets, chapter 3.3, pp. 799-860, in Vol B of {\it
Handbook of Convex Geometry}, ed.\ P. Gruber and J. Wills, North
Holland, Amsterdam, 1993.

\item{[FKK]} G. Fejes T\'oth, G. Kuperberg and W. Kuperberg, Highly
saturated packings and reduced coverings, {\it Monatsh.  Math.}  
125(1998) 127-145.

\item{[GrW]} P. Gruber and J. Wills, eds., {\it Handbook of Convex Geometry}, 
North Holland, Amsterdam, 1993.

\item{[GrS]}\ B. Gr\"unbaum  and G.C. Shephard, {\it Tilings and Patterns},
Freeman, New York, 1986.

\item{[Kup]} G. Kuperberg, Notions of denseness, {\it Geom. Topol.} 4(2000)
274-292.

\item{[MaM]} G.A. Margulis and S. Mozes, Aperiodic tilings of the
hyperbolic plane by convex polygons, {\it Israel J. Math.}
107(1998) 319-332.

\item{[Nev]} A. Nevo, Pointwise ergodic theorems for radial averages
on simple Lie groups I, {\it Duke Math. J.} 76(1994) 113-140.

\item{[NeS]} A. Nevo and E. Stein, Analogs of Weiner's ergodic
theorems for semisimple groups I, {\it Annals of Math.} 145(1997)
565-595.

\item{[Pen]} R. Penrose, Pentaplexity - a class of non-periodic tilings of the
plane, {\it Eureka} 39(1978) 16-32. (Reproduced in {\it Math. Intell.} 2(1979/80) 32-37.)

\item{[Rad]} C. Radin, {\it Miles of Tiles}, Student Mathematical Library,
Vol 1,
Amer. Math. Soc., Providence, 1999.

\item{[RaW]} C. Radin and M. Wolff, Space tilings and local isomorphism,
{\it Geometriae Dedicata} 42 (1992) 355-360.

\item{[Rog]} C.A. Rogers, {\it Packing and Covering}, University Press,
Cambridge, 1964.

\vfill \eject
\hbox{}
\vs.5
\hs0 \epsfig 1\hsize; 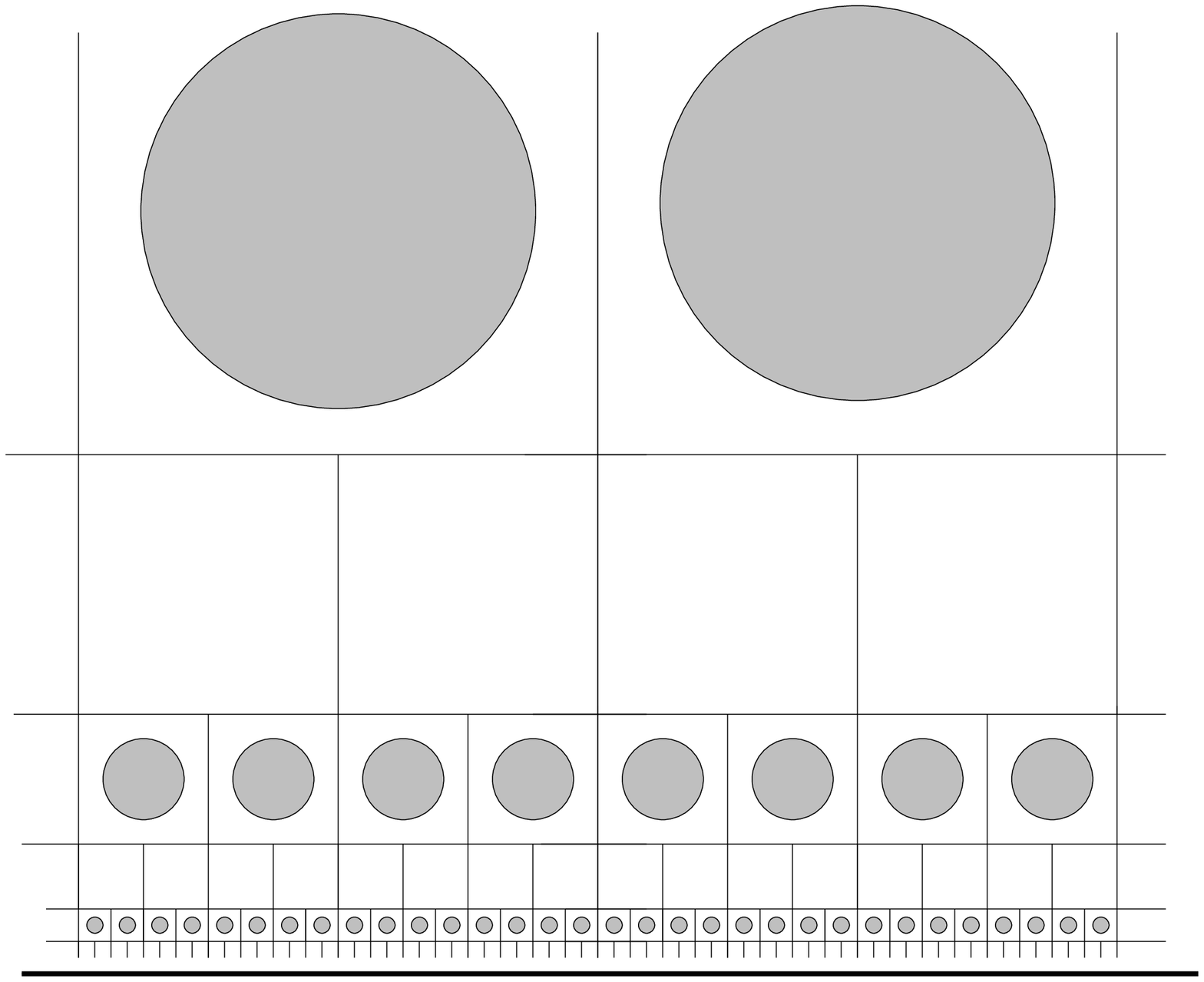 ;
\vs.3 \centerline{Figure 1. Boroczky's packing of disks}
\vfill \eject
\hbox{}
\vs.5
\hs0 \epsfig 1\hsize; 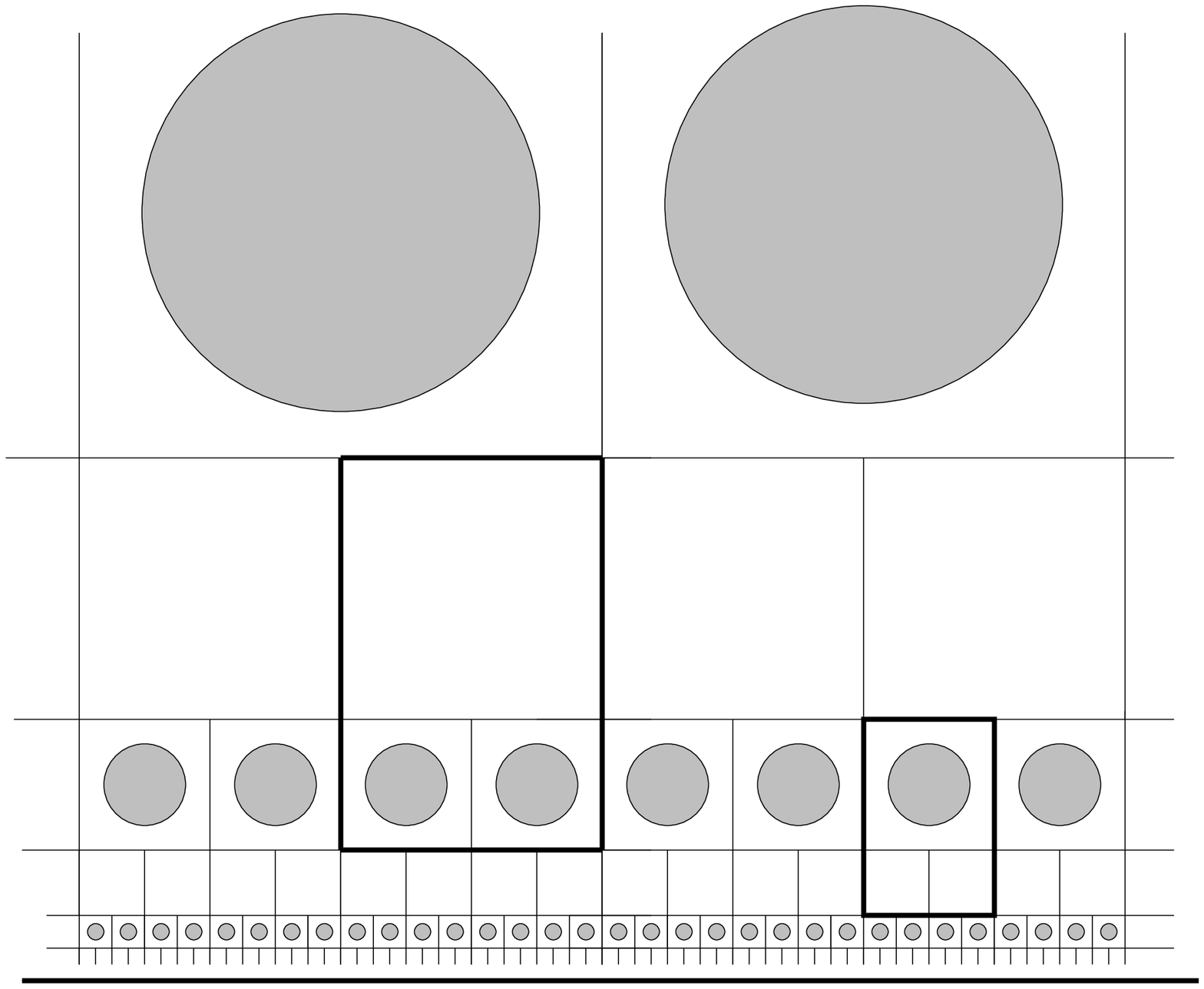 ;
\vs.3 \centerline{Figure 2. Boroczky's packing with two tiles in dark
outline}
\vfill \eject
\hbox{}
\vs-.5
\hs0 \epsfig .25\hsize; 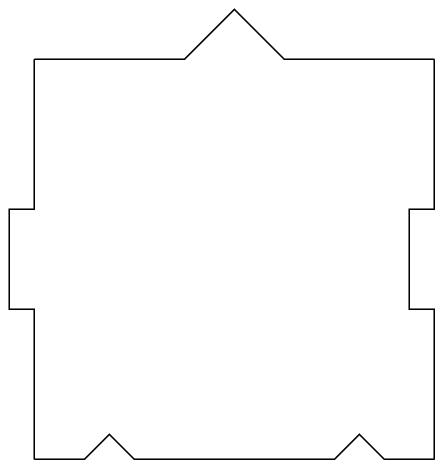 ;
\vs.2 \centerline{Figure 3. A tile}
\vs.4
\hs0 \epsfig 1\hsize; 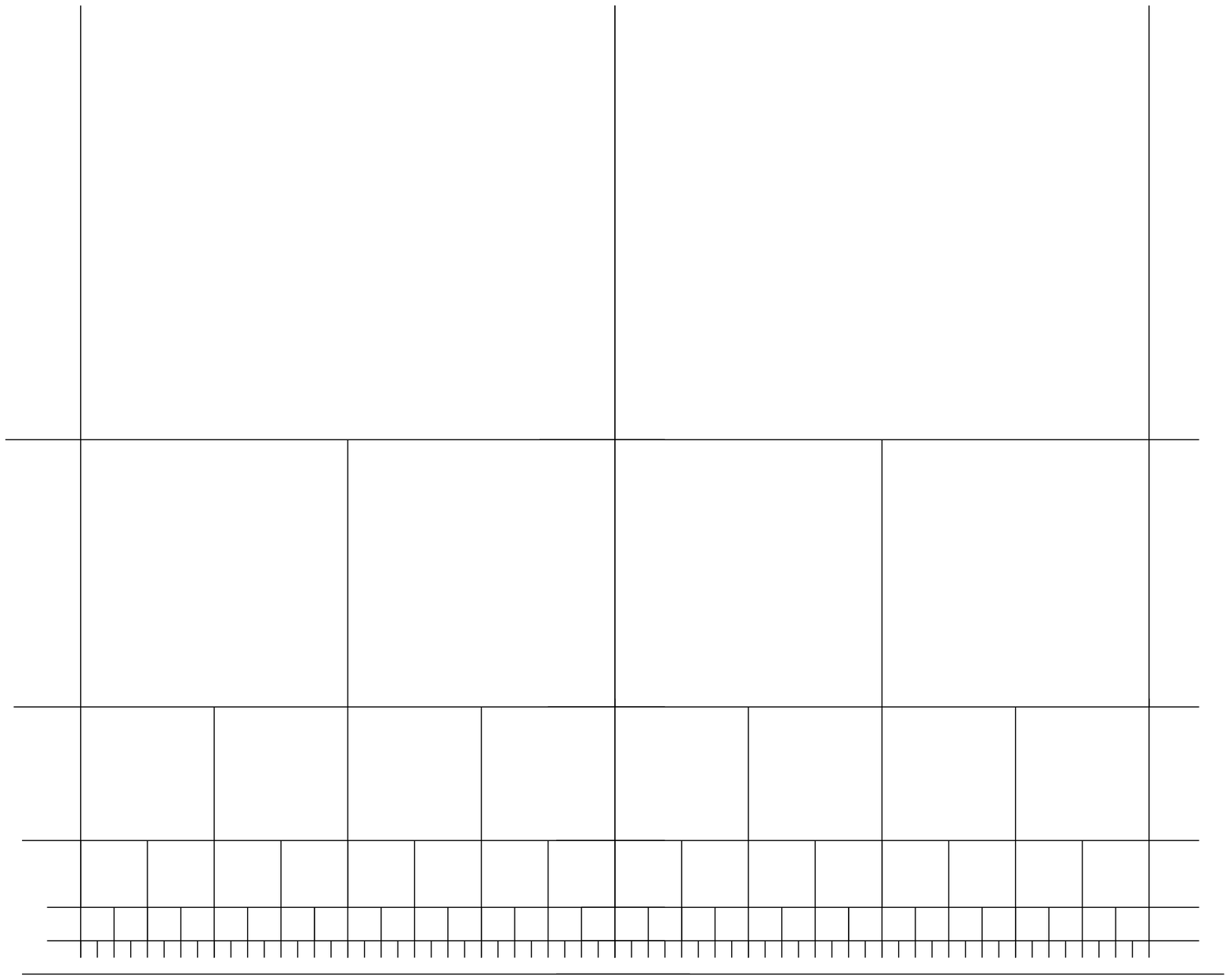 ;
\vs.3 \centerline{Figure 4. A tiling}
\vfill \eject
\hbox{}
\vs0
\hs0
\epsfig .25\hsize; 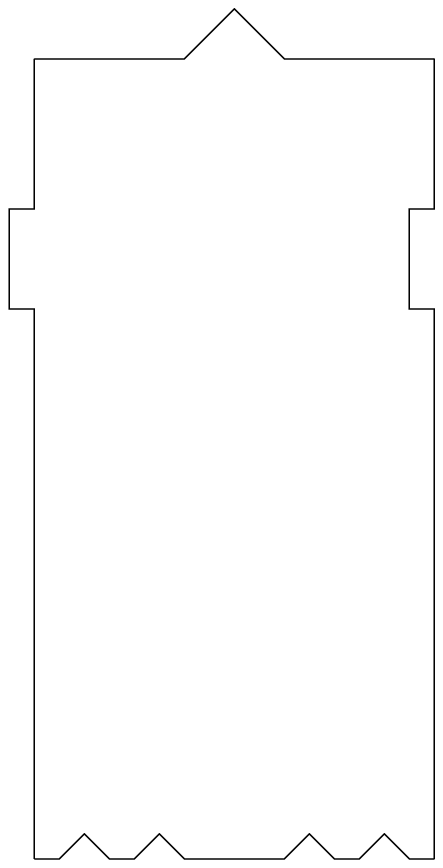 ;
\vs.3 \centerline{Figure 5. A tile}
\vs.5
\hs0 \epsfig .25\hsize; 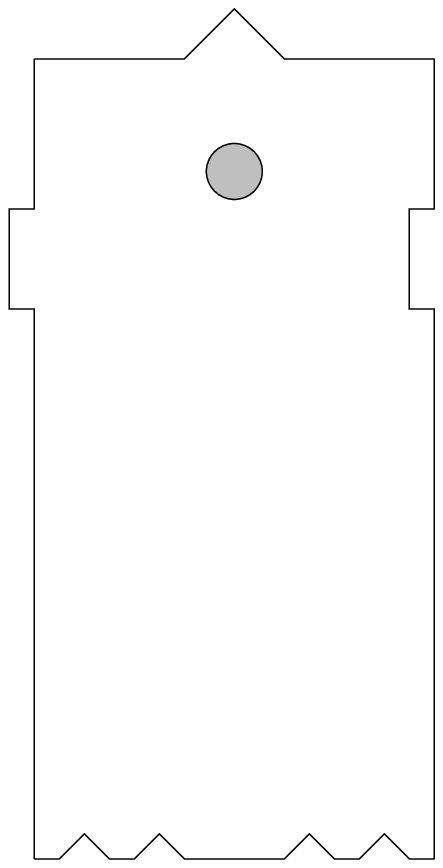 ;
\vs.3 \centerline{Figure 6. A body}
\vfill \eject
\hbox{}
\vs-.2
\hs0 \epsfig .8\hsize; 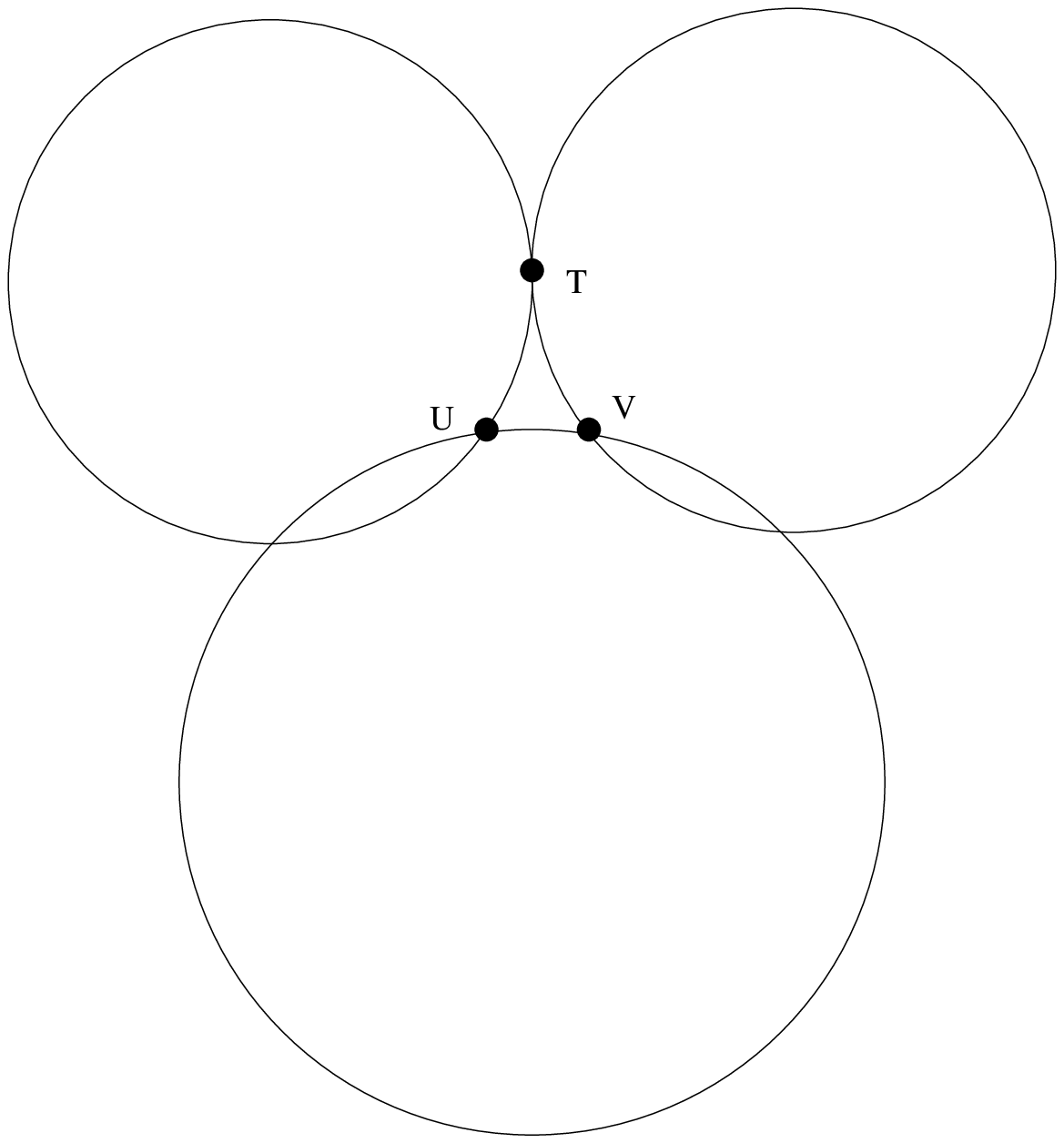 ;
\vs.2 \centerline{Figure 7. Uncovered regions}
\vfill \eject
\nd
\hbox{}
\vs-.2
\hs0 \epsfig .4\hsize; 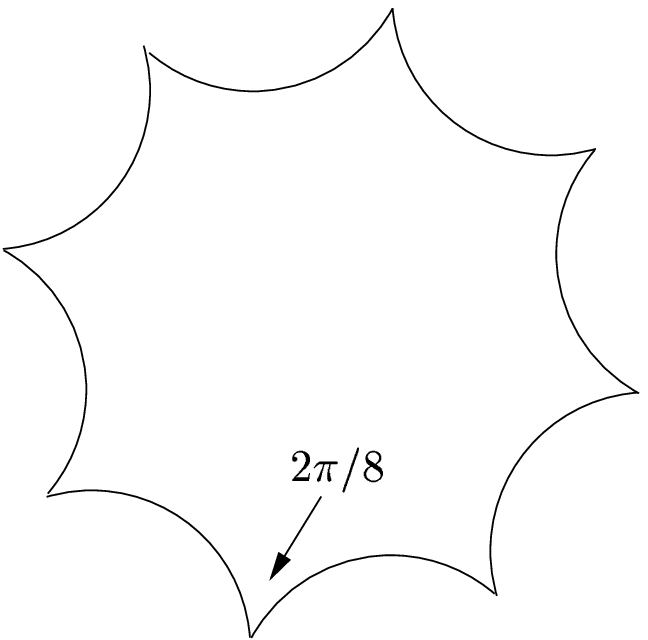 ;
\vs.2 \centerline{Figure 8. $\beta_1$}
\vs.5
\hs0 \epsfig .4\hsize; 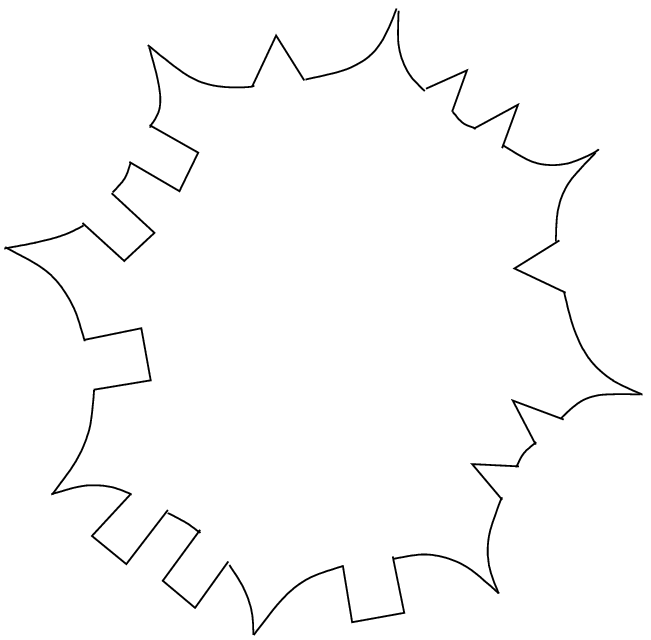 ;
\vs.2 \centerline{Figure 9. $\beta'_2$}
\vfill \eject
\hbox{}
\vs-.2
\hs0 \epsfig .4\hsize; 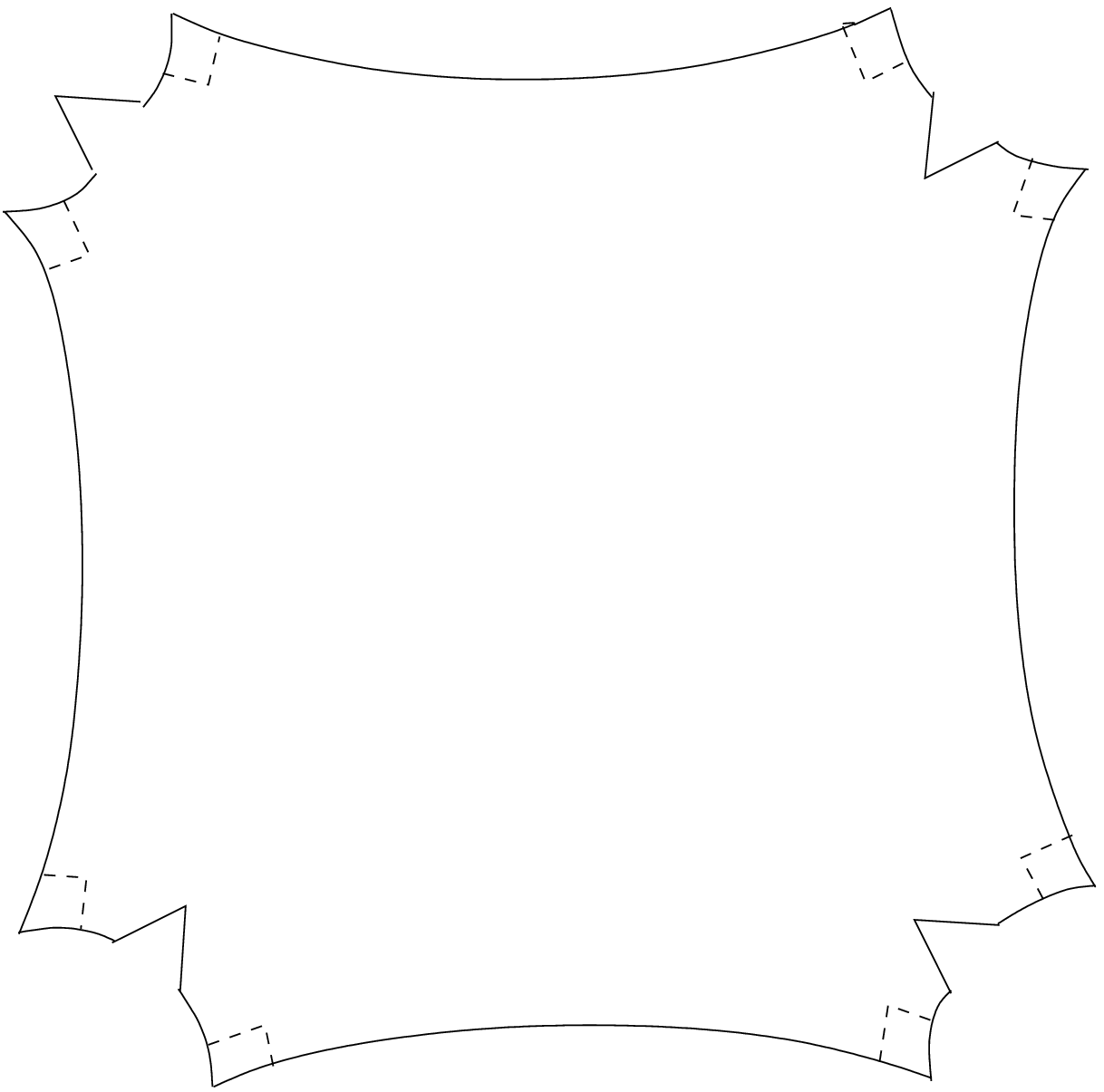 ;
\vs.2 \centerline{Figure 10. $\tau$}
\vs.5
\hbox{}
\vs-.2
\hs0 \epsfig .4\hsize; 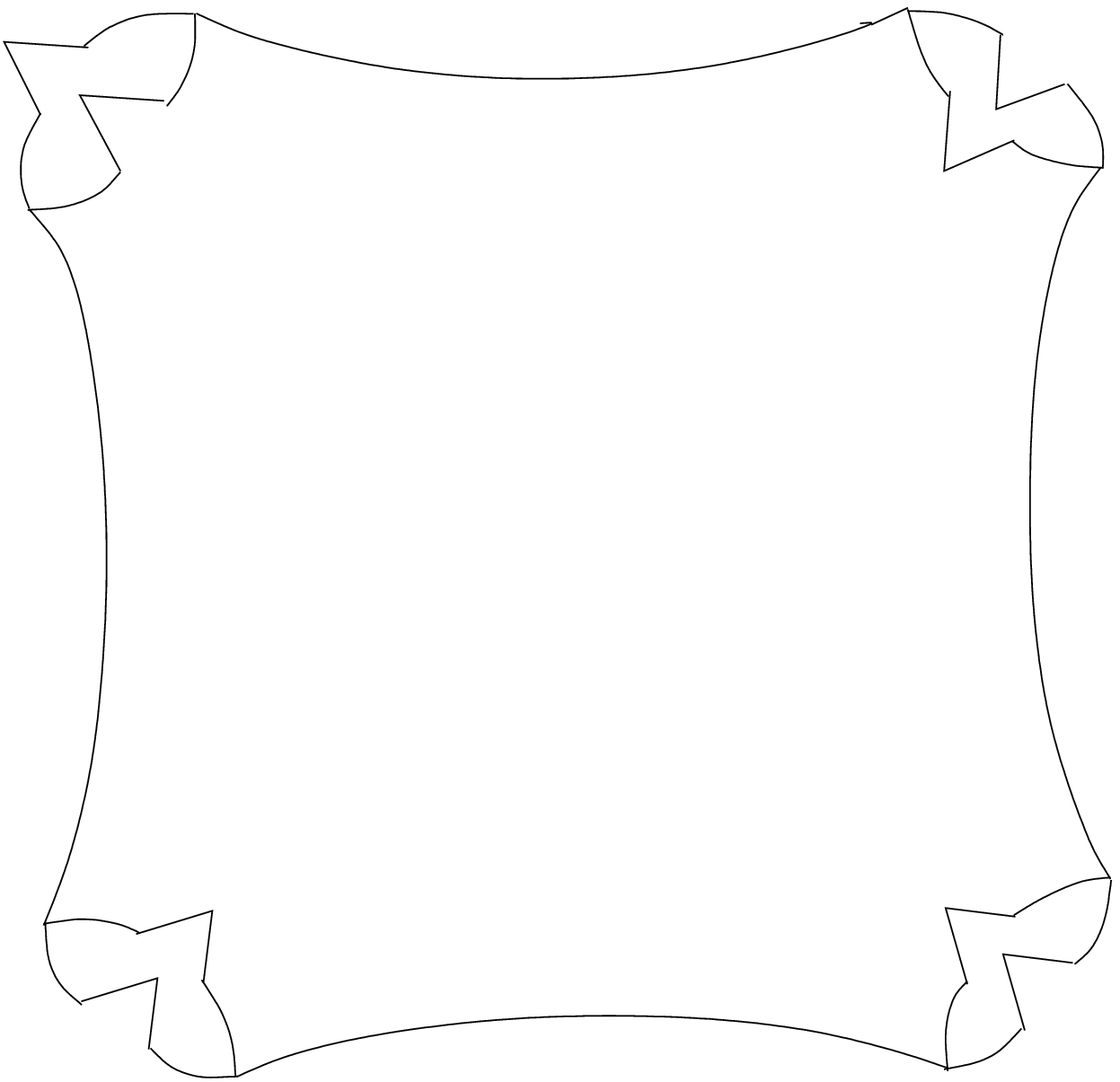 ;
\vs.2 \centerline{Figure 11. $X$}
\vfill \eject
\end